\documentstyle[12pt]{article} \textwidth 6.1 in \textheight 9.5 in

\date{}
\title{\bf  Uniqueness of singular solution of
semilinear elliptic equation
 \thanks{ Supported in part by the National Science Funds of China
 (10901047)} \\
 }
\author{Baishun, Lai and  Qing, Luo   \\
\small {\it Department of Mathematics, Henan University, Kaifeng, 475001, PR China }\\
\small {\it E-mail: laibaishun@henu.edu.cn }\\
} \baselineskip 0.2in
\date{}
\begin{document}
\maketitle
\begin{center}
\begin{minipage}{130mm} {\small {\bf Abstract}

\ \  \ In this paper, we study asymptotic behavior of solution near
0 for a class of elliptic problem. The uniqueness of singular
solution is established.\vskip0.1in

 {\it Keywords:}\
Nonhomogeneous semilinear elliptic equation; Positive solutions;
  Asymptotic behavior; Singular solutions \ \ }
\end{minipage}
\end{center}
\baselineskip 0.2in \vskip 0.2in {\bf 1. Introduction} %
\setcounter{section}{1} \setcounter{equation}{0}\vskip 0.1in In this
paper, we  study the elliptic equation
\begin{equation}
\Delta u+K(|x|)u^{p}+\mu f(|x|)=0,
\end{equation}
where $n\geq 3, \Delta=\Sigma_{1}^{n}\frac{\partial^{2}}{\partial
x_{i}^{2}}$ is the Laplace operator, $p>1$ is a constant, $\mu\geq0$
is a parameter, and $f$ and $K$ are given locally H$\ddot{o}$lder
continuous function in $R^{n}\setminus\{0\}$, so that the solutions
of (1.1) are classical on $0<|x|<\infty$. However, at $x=0$, when
$K$ is "bad", usually we cannot expect the solutions to be
differentiable, or even continuous owing to the singularity of $K$
at $x=0$. Let $u$ be a solution of (1.1), the singular point $x=0$
of $u$ is called a removable singular point if
$u(0)\equiv\lim_{x\to0}u(x)$ exists, otherwise $x=0$ is called a
nonremovable singular point. It is shown by Ni and Yotsutani [13]
that when $x=0$ is a removable singular point of a solution of
(1.1), the existence of the derivatives of the solution depends on
the ''blow up" rate of $K$ at $x=0$ [13, Proposition 4.4].

Let $u\in C^{2}(R^{n}\setminus 0)$ be a solution of (1.1). If $x=0$
is a removable singular point of $u$, then $u$ is said to be a
regular solution of (1.1), otherwise $u$ is said to be a singular
solution.

The purpose of this paper is to study the asymptotic behavior of
singular positive solutions and to obtain the uniqueness of singular
positive solutions of (1.1) which has diverse physical and
geometrical backgrounds. In particular,  in the case $K=1$ and
$p=2$, (1.1) arises naturally in establishing occupation time limit
theorems for super-Brownian motions which requires analyzing
cumulant generating functions satisfying some integral equations
equivalent to the parabolic counterparts of (1.1). There are many
works devoted to the studying the existence, monotonicity and
asymptotic expansion at infinity of positive solutions of the
equation (1.1). We refer the interested readers to [2,3,4,6,7, 9-12]
and the references therein.

 In this paper, we consider positive radial solutions of (1.1) with
  radial functions $K,f$. The radial version of (1.1) is of the form
$$u''+\frac{n-1}{r}u'+K(r)u^{p}+\mu f=0.$$

For the same reasons, the regular solutions that have finite limits
at $r=0$, are particularly interesting , which lead us to consider
the initial value problem
$$
\left\{
\begin{array}{ll}u''+\frac{n-1}{r}u'+K(r)u^{p}+\mu f=0,\\
u(0)=\alpha.
\end{array}
\right. \eqno (1.2)
$$
We use $u_{\alpha}=u(r,\alpha)$ to denote the unique solution of
(1.2).
 \vskip0.1in

First, we introduce the following notations, which will be used
throughout this paper:
\[
m\equiv \frac{l+2}{p-1},\ \ \ \ \ L\equiv \lbrack
m(n-2-m)]^{\frac{1}{p-1}},
\]%

\begin{eqnarray*}
p_{c}=\left\{
\begin{array}{ll}
\frac{(n-2)^{2}-2(l+2)(n+l)+2(l+2)\sqrt{(n+l)^{2}-(n-2)^{2}}}{(n-2)(n-10-4l)},
& n>10+4l,  \\
\infty,& 3\leq n\leq 10+4l.
\end{array}
\right.
\end{eqnarray*}%
$$
\lambda _{2}=\lambda _{2}(n,p,l)=\frac{(n-2-2m)+\sqrt{
(n-2-2m)^{2}-4(l+2)(n-2-m)}}{2}.
$$%

The hypotheses of $K(|x|)$ are often divided into two cases: the
fast decay case and the slow decay case. In this paper, we will
focus on the slow decay case, i.e. $K(r)\geq cr^{l}$, for some
$l>-2$ and $r$ large. First, let us introduce a collection of
hypotheses on $K(|x|)$ and $f$:

$(K.1)$\ \ $K(|x|)=k_{\infty }|x|^{l}+O(|x|^{-d_{1}})$ at
$|x|\rightarrow \infty $ for some constants $k_{\infty }>0$ and
$d_{1}>n-\lambda _{2}-m(p+1)$.

$(K.2)$ \ \ $\lim_{x\to0}|x|^{-l}K(|x|)=k_{0}>0.$

$(K.3)$\ \ $K(r)$ is locally Lipschitz continuous and $\frac{d}{dr}
(r^{-l}K(r))\leq0$ for a.e. $r>0$.

$(f.1)$\ \ $\lim_{|x|\to 0} |x|^{d}f(|x|)=b\geq0$, where $0<d\leq
m+2$.

$(f.2)$\ \ $f(x)=O(|x|^{-q})$ near $\infty $ for some $q>n-m-\lambda
_{2}$.

$(f'.2)$ $\int_{0}(r^{d}f)_{r}^{+}dr<\infty.$

$(K'.2)$ $\int_{0}(r^{-l}K(r))_{r}^{+}dr< \infty.$

$(f'.3)$ $\int_{0}(r^{d}f)_{r}^{-}dr<\infty.$

$(K'.3)$ $\int_{0}(r^{-l}K(r))_{r}^{-}dr< \infty,$ where
$k^{\pm}=\max\{\pm k,0\}, r=|x|$. \vskip0.1in

Our main result is  as follows:\vskip0.1in

{\bf Theorem 1.}  Suppose that $K(r)$ satisfies $(K.1)-(K.3),(K'.2)$
$f$ satisfies $(f.1)$ and $(f.2),(f'.2)$, $p>p_{c},0<d<2$. Then
(1.2) has one and only one singular solution $U(r)$, Furthermore,
for any regular solution $u(r)$, the following holds
$$u(r)<U(r)\leq\frac{L}{(r^{-l}K(r))^{\frac{1}{p-1}}r^{m}}.$$\vskip0.1in

This paper is organized as follows. In Section 2, the asymptotic
behavior of singular positive solution of (1.2) near 0 is studied.
Finally, the uniqueness result about  the singular positive solution
is established.\vskip0.1in

{ \bf 2  Asymptotic behavior of singular solution near 0}\setcounter{section}{2} %
\setcounter{equation}{0}\vskip0.1in In this section, we obtain the
asymptotic behavior of the positive radial solution of (1.2) near 0.
First, we obtain the prior estimates of the positive solution of
(1.2) near 0.

{ \bf Lemma 2.1.} Let $p>1$ and $u$ be a positive solution of (1.2)
in $B_{r}\setminus\{0\}$. Then, there exists a positive constant $R$
such that for $0<r<R$, $u(r)\leq C r^{-\frac{2+l}{p-1}}$ for some
constant $C$.\vskip0.1in

{\bf Proof.} From (1.2), we have
$$(r^{n-1}u'(r))'+r^{n-1}K(r)u^{p}\leq0,$$ then, there exists a small
positive constant $r_{k}<r$ such that \begin{eqnarray*}
r^{n-1}u'(r)&=&r^{n-1}_{k}u'(r_{k})-\int_{r_{k}}^{r}r^{n-1}K(r)u^{p}dr\\
&\leq&-\int_{r_{k}}^{r}r^{n-1}K(r)u^{p}dr,
\end{eqnarray*}
since $u'(r_{k})<0$ near 0. From that,
$$\frac{u'(r)}{u^{p}(r)}\leq-\frac{C}{r^{n-1}}\int_{r_{k}}^{r}s^{n-1}K(s)ds\leq
Cr^{l+1} \eqno(2.1).$$ Integrating (2.1) over $(0,r)$, we have
$u(r)\leq Cr^{-\frac{2+l}{p-1}}$, and the proof is
completed.\vskip0.2in

Now, we verify the following asymptotic behavior of $u$ near 0 by
using the Li's energy method in [9]\vskip0.2in

{ \bf Theorem 2.2.} Let $u$ be a positive solution of (1.2) near 0.
Assume that $u(r)=O(r^{-m})$ at 0 and $ K, f$ satisfy

(i) $(f.1), (K.1)$ and $(f'.2), (K'.2)$ if $p>\frac{n+2+2l}{n-2}$ or

(ii) $(f.1),(K.1)$ and $(f'.3), (K'.3)$ if
$\frac{n+l}{n-2}<p<\frac{n+2+2l}{n-2}$ with $d=m+2, b\geq0$.

 Then, $b\leq\max_{z\in R^{+}}\{L^{p-1}z-k_{0}z^{p}\}$ and $\lim_{r\to0}r^{m}u(r)=z_{1}\
\mbox{or}\ z_{2}$, where $z_{1}$ and $z_{2}$, $z_{1}\leq z_{2}$ are
two roots of the equation $k_{0}z^{p}-L^{p-1}z+b=0$.\vskip0.1in

{\bf Proof.}  Denote $v(t):=r^{m}u, t=\log r,$ then
$$v''(t)+av'(t)-L^{p-1}v+k(t)v^{p}+g(t)=0, \eqno(2.2)$$ where
$a=n-2-2m, k(t)=e^{-lt}K(e^{t})$ and $g(t):=e^{(m+2)t}f(e^{t})$.
Suppose that $$0\leq
\alpha=\lim\inf_{t\to-\infty}v(t)<\lim\sup_{t\to-\infty}v(t)=\beta<\infty.$$
Then, there exist two sequences $\{\eta_{i}\}$ and $\{\xi_{i}\}$
going to $-\infty$ as $i\to-\infty$ such that $\{\eta_{i}\}$ and
$\{\xi_{i}\}$ are local minima and local maxima of $v$,
respectively, satisfying $\eta_{i}<\xi_{i}<\eta_{i+1}, i=1.2,....$
Define an energy function
$$E(t):=\frac{1}{2}(v')^{2}-\frac{L^{p-1}}{2}v^{2}+\frac{1}{p+1}k(t)v^{p+1}+bv. \eqno(2.3)$$
Now, multiplying (2.2) by $v'(t)$ and integrating by parts over
$[t,T]$ ($T$ is a constant) we obtain \[
E(t)+a\int_{t}^{T}(v')^{2}ds+
[g(t)-b]v(t)=C(T)+\int_{t}^{T}g'ds+\frac{1}{p+1}\int_{t}^{T}v^{p+1}k'ds.
\ \ \ \ \ \ \ \ \ \ \ \ \ \ \ \ \ \ \ \ \ \ \ \ \ (2.4)
\]
 Assume that $f, K$ satisfy
$(f.1), (K.1)$ and $(f'.2), (K'.2)$ and $p>\frac{n+2+2l}{n-2}$, from
(2.4), we have \begin{eqnarray*} &&E(t)+a\int_{t}^{T}(v')^{2}ds+
[g(t)-b]v(t)+\frac{1}{p+1}\int_{t}^{T}v^{p+1}[k']^{-}ds+\int_{t}^{T}[g']^{-}ds\\
&&=C(T)+\int_{t}^{T}[g']^{+}ds+\frac{1}{p+1}\int_{t}^{T}v^{p+1}[k']^{+}ds.
\end{eqnarray*}
Since $v$ is bounded, $E(\eta_{i})$ is bounded independently of $i$
and $[k']^{+}, [g']^{+}\in L^{1}(T,\infty)$, we conclude that
$$\int_{T}^{\eta_{i}}V_{s}^{2}ds,\ \ \int_{T}^{\eta_{i}}
[k']^{-}V^{p+1}ds\ \ \mbox{and} \ \int_{T}^{\eta_{i}}[g']^{-}ds$$
are bounded independent of $i$ which implies that
$$\int_{-\infty}^{T}(v')^{2}ds<\infty. \eqno(2.5)$$
Similarly, if $f, K$ satisfy $(f.1), (K.1)$ and $(f'.3), (K'.3)$, we
also have $V_{s}\in L^{2}(-\infty,T)$. It follows from (2.4) and
(2.5) that $E=\lim_{t\to-\infty}E(t)$ exists, which in turn from
(2.3) implies $v'(t)$ is bounded. Then, from (2.2), $v''(t)$ is
bounded also and from (2.5), $v'(t)\to 0$ as $t\to -\infty$.

Let $h(v)=-\frac{L^{p-1}}{2}v^{2}+\frac{k_{0}}{p+1}v^{p+1}+bv$.
Since
$$\lim_{i\to\infty}E(\eta_{i})=h(\alpha)=E=h(\beta)=\lim_{i\to\infty}E(\xi_{i}),$$
we  choose $\alpha<\gamma<\beta$ and $t_{i}\in(\eta_{i},\xi_{i})$
such that $v(t_{i})=\gamma, \frac{dh}{dv}(\gamma)=0$ and
$h(\gamma)\neq E$. However
$E=\lim_{i\to\infty}E(t_{i})=\lim_{i\to\infty}(\frac{1}{2}v'(t_{i})^{2}+h(\gamma))=h(\gamma)$,
 a contradiction. Therefore,
$v_{\infty}=\lim_{t\to-\infty}v(t)$ exists. For given
$\tilde{\xi}>0$, there exists a sequence $\{s_{i}\}$ converging to
$-\infty$ such that $|v'(s_{i})|\leq \tilde{\xi}$ for $ i=1,2,....$
Since $E(s_{i})$ is bounded, we obtain (2.5) from (2.4). From (2.4)
and (2.5), $\lim_{t\to-\infty}E(t)$ exists. Thus, (2.3) implies
$\lim_{t\to-\infty}v'(t)=0$. Then, by (2.2), $\lim_{t\to\infty}
v''(t)$ exists and must be 0. Therefore, we conclude from (2.2) that
$b\leq\max_{z\in R^{+}}\{L^{p-1}z-k_{0}z^{p}\}$  and
$v_{\infty}=z_{1}\ \mbox{or}\ z_{2}$. The proof is completed.
\vskip0.2in

{\bf Corollary 2.3.} Let $u$ be a positive solution of (1.2) near 0,
and $f, K$ satisfy the same condition as in Theorem 2.2 except that
$d=m+2$ is replaced by $0<d<m+2$. Then
$\lim_{r\to0}r^{m}u(r)=\frac{L^{p-1}}{k_{0}}$ or 0.\vskip0.1in In
case $0<d<m+2$, then $\lim_{t\to\infty}
g(t)=\lim_{t\to-\infty}e^{(m+2)t}f(e^{t})=0$, similar to that of
Theorem 2.2, we can immediately obtain
$\lim_{r\to0}r^{m}u(r)=\frac{L^{p-1}}{k_{0}}$ or 0. The detail proof
is omitted here

 \vskip 0.2in

{\bf Theorem 2.4.} \ Let $p>\frac{n+l}{n-2}$. Suppose that
$f(r)=O(r^{-d})$ at 0 with $d<m+2$. Then, any positive solution of
(1.2) satisfying $\lim_{r\to 0}r^{m}u(r)=0$ has the asymptotic
behavior at 0 such that
\begin{eqnarray*}
u(r)=\left\{ \begin{array}{lll} O(r^{2-d}) & \mbox{if}\ 2<d<m+2,\\
O(|\log r|) & \mbox{if}\ d=2,\\
O(1) & \mbox{if}\ d<2.
 \end{array}
 \right.
\end{eqnarray*}
\vskip 0.2in

{\bf Proof.} First, we claim that there exists a constant $\beta>0$
such that  $u(r)=O(r^{-m+\beta})$ near 0. Set $v(r)=r^{m}u(r)$ for
$r>0$ and $L_{\varepsilon}v \equiv \Delta v-\frac{2mv^{\prime
}}{r}-m(n-2-m-\varepsilon )\frac{v}{r^{2}}+\mu r^{m}f.$ Then, $v(r)$
satisfies
$$L_{\varepsilon}v-m\varepsilon\frac{v}{r^{2}}+\frac{v^{p}}{r^{2}}K(r)r^{-l}=0$$
and for any $\varepsilon>0$, there exists $R_{\varepsilon}>0$ such
that $L_{\varepsilon}v\geq0$ for $0<r\leq R_{\varepsilon}$.

On the other hand, for $0<\varepsilon <n-2-m$, let $\varphi
_{\varepsilon
}(x)=|x|^{\beta _{\varepsilon }}$ we have%
\[
L_{\varepsilon }\varphi _{\varepsilon }=[\beta (\beta
-1)+(n-1-2m)\beta -m(n-2-m-\varepsilon )]|x|^{\beta _{\varepsilon
}-2}+\mu r^{m}f
\]%
in $R^{n}\setminus \{0\}.$ Choosing $\beta _{\varepsilon }>0$
sufficiently small such that
\[
\beta _{\varepsilon }(\beta _{\varepsilon }-1)+(n-1-2m)\beta
_{\varepsilon }-m(n-2-m-\varepsilon )\leq 0,
\]%
and
\[
\frac{r^{m}f}{r^{\beta _{\varepsilon }-2}}\rightarrow 0 \ \ \ \ \
\mbox{as}\ \ \ \ \ r \to0 .
\]%
So there exists an $R_{\varepsilon }^{\prime }>0$ such that
\[
L_{\varepsilon }\varphi _{\varepsilon }\leq 0\ \ \ \ \ \ \mbox{in}\
\ \ \ \ \ \ \ 0<r<R_{\varepsilon }^{\prime }.
\]%
Setting $R_{\varepsilon }^{\prime \prime }=\min \{R_{\varepsilon
}^{\prime },R_{\varepsilon }\}$, $C_{\varepsilon }=v(R_{\varepsilon
}^{\prime \prime })(R_{\varepsilon }^{\prime \prime })^{-\beta
_{\varepsilon }},$ we see that
\begin{eqnarray*}
\left\{
\begin{array}{lll}
L_{\varepsilon }(v-C_{\varepsilon }\varphi _{\varepsilon })\geq 0 &
\ \ \ \ \ \ \mbox{for}\ \ \ \ 0<r\leq R_{\varepsilon }^{\prime \prime }, \\
v-C_{\varepsilon }\varphi _{\varepsilon }=0 & \ \ \ \ \ \ \mbox{at}\
\ \ \
R_{\varepsilon }^{\prime \prime }, \\
v(r)-C_{\varepsilon }\varphi _{\varepsilon }(r)\rightarrow 0 & \ \ \ \ \ \ %
\mbox{as}\ \ \ \ \ r\to 0 , &
\end{array}
\right.
\end{eqnarray*}
since $\beta _{\varepsilon }>0$. Observing that the coefficient of the term $%
v$ in $L_{\varepsilon }$ is negative, we conclude by the maximum
principle that $v-C_{\varepsilon }\varphi _{\varepsilon }\leq 0$ for
$0<r\leq R_{\varepsilon }^{\prime \prime }$, i.e. $v(r)\leq
C_{\varepsilon }r^{\beta _{\varepsilon }}$ near 0. This guarantees
that $u(r)\leq C_{\varepsilon }r^{-m+\beta _{\varepsilon }}$ near 0.
Since, there exists a sequence ${t_{i}}$ going to $-\infty$  such
that $v'(t_{i})\to 0$  with $r_{i}=e^{t_{i}}$, then,
$r_{i}^{\frac{p+1+l}{p-1}}u'(r_{i})\to0$ as $r_{i}\to0$  and thus $$
\lim_{r_{i}\to0}r_{i}^{n-1}u'(r_{i})=0. \eqno(2.6)$$ By  (1.2) and
(2.6), we observe that
$$u'(r)=-\frac{1}{r^{n-1}}\int_{0}^{r}(K(r)u^{p}+\mu f)s^{n-1}ds
\eqno(2.7)$$ for $0<r\leq R$. Integrating (2.7) over $[r,R]$, and
changing the order of integration, we have
$$\begin{array}{lllllllll} u(r)&=&
\int_{r}^{R}\frac{1}{t^{n-1}}[\int_{0}^{t}(K(s)u^{p}+\mu
f)s^{n-1}ds]dt+u(R)\\\\
&\leq& \frac{r^{2-n}}{n-2}\int_{0}^{r}(K(s)u^{p}+\mu f)s^{n-1}ds+
\frac{1}{n-2}\int_{r}^{R}(K(s)u^{p}+\mu f)sds+u(R).
\end{array}\eqno(2.8)$$

Then, for small  $r>0$, \begin{eqnarray*} u(r)\leq \left\{
\begin{array}{lll} C(1+r^{-m+p\beta}+r^{2-d} & \mbox{if}\
p\beta\neq m,d>2,\\
C(1+|\log r|+r^{2-d}) & \mbox{if}\ p\beta=m,d>2.
\end{array}
\right.
\end{eqnarray*}
If $m+2-p\beta\leq d$, then $u(r)=O(r^{2-d})$ near 0. Otherwise,
repeating the above arguments with $\beta$ replace by $p\beta$ leads
to $m+2-p^{2}\beta>d$ and in turn $m+2-p^{k}\beta>d$ for any
positive integer $k$, which is absurd. Therefore, we have
$u(r)=O(r^{2-d})$ near 0.

If $d=2$, then it follows (2.8) that for small $r>0$,
\begin{eqnarray*} u(r)\leq \left\{
\begin{array}{lll} C(1+r^{-m+p\beta}+|\log r|) & \mbox{if}\
p\beta\neq m, \\
C(1+|\log r|) & \mbox{if}\ p\beta= m.
\end{array}
\right.
\end{eqnarray*}
The case $p\beta\neq m$ can be dealt as same in the above.
Therefore, $u(r)=O(\log r)$ near 0. On the other hand, if $d<2$,
then it follows (2.8) that for small $r>0$,
\begin{eqnarray*} u(r)\leq \left\{
\begin{array}{lll} C(1+r^{-m+p\beta}) & \mbox{if}\
p\beta\neq m, \\
C(1+|\log r|) & \mbox{if}\ p\beta= m.
\end{array}
\right.
\end{eqnarray*}
If $-m+p\beta>0$, then $u(r)=O(1)$ near 0. Otherwise, repeating the
above arguments with $\beta$ replace by $p\beta$ leads to
$-m+p^{2}\beta<0$ and in turn $-m+p^{k}\beta<0$ for any positive
integer $k$, which is absurd. Therefore, we have $u(r)=O(1)$ near 0.
The proof is completed.

\baselineskip 0.2in \vskip 0.2in  %
{\bf 3  Uniqueness of singular solution }\setcounter{section}{5} %
\setcounter{equation}{0}\vskip0.1in

In order to prove the Theorem 1, we need following Lemma\vskip0.1in

{\bf Lemma 3.2.} Suppose that $K(r),f(r)$ satisfy the same condition
as in the Theorem 2.2 except that $d=m+2$ is replaced by $d<m+2$.
Let $u(r)$ be a positive solution of (1.2), then,
$$\lim_{r\to0^{+}}r\frac{d}{dr}(r^{m}u(r))=0.$$

{\bf Proof.} Let $v(t)=r^{m}u(r), r=e^{t}$, then
$$v''+b_{0}v'-L^{p-1}v+k(t)v^{p}+\mu f(e^{t})e^{(m+2)t}=0$$ or
$$(e^{b_{0}t}v'(t))'+e^{b_{0}t}(k(t)v^{p}-L^{p-1}v+\mu
f(e^{t})e^{(m+2)t})=0.$$ Integrating from $T$ to $t$ for $T<t$, we
have
$$e^{b_{0}t}v'(t)-e^{b_{0}T}v'(T)+\int_{T}^{t}e^{b_{0}s}(k(s)v^{p}-L^{p-1}v+\mu
f(e^{s})e^{(m+2)s})=0.$$ By Corollary 2.3, we know that
$\lim_{t\to-\infty}(k(t)v^{p}-L^{p-1}v+\mu f(e^{t})e^{(m+2)t})=0$.
Given $\varepsilon>0$, there exists $t_{\varepsilon}$ such that
$|k(t)v^{p}-L^{p-1}v+\mu f(e^{t})e^{(m+2)t}|<\varepsilon$, when
$t<t_{\varepsilon}$, and $$|v'(t)|\leq
|e^{b_{0}(T-t)}v'(T)|+\frac{\varepsilon}{b_{0}}(-e^{-b_{0}(t-T)}+1).$$
Letting $T$ (if necessary, take a subsequence) go to $-\infty$, we
have $|v'(t)|\leq\frac{\varepsilon}{b_{0}}$, if $t<t_{\varepsilon}$,
since $\varepsilon$ can be arbitrary small, we conclude
$\lim_{t\to-\infty}v'(t)=0$, or equivalently
$\lim_{r\to0^{+}}r\frac{d}{dr}(r^{m}u(r))=0.$\vskip0.1in

{\bf Lemma 3.3.}$^{[11]}$ \ Suppose $f(t)$ and $g(t)$ are continues
functions, $\lim_{t\to+\infty}f(t)=b>0$,
$\lim_{t\to+\infty}g(t)=c>0$. Let $y(t)$ be a solution of
$$y''-f(t)y'+g(t)y=0.$$ Then $y(t)$ is unbounded as
$t\to+\infty$.\vskip0.1in

{\bf Proof of Theorem 1.}  First, we prove (1.2) has only one
singular solution. Suppose $u_{1}(r)$ and $u_{2}(r)$ are two
different singular solutions. As we did in Lemma 3.2, let
$v_{i}(t)=r^{m}u_{i}(r), i=1,2, r=e^{t}$, and
$h(t)=\frac{v_{2}(t)}{v_{1}(t)},$ then we have $$
h''+(b_{0}+\frac{2v_{1}'}{v_{1}})h'+k(t)v_{1}^{p-1}(h^{p}-h)-\mu
f(e^{t})e^{(m+2)t}(\frac{h-1}{v_{1}})=0.$$ Let $Q(t)=h(-t)-1,
f(t)=b_{0}+\frac{2v_{1}'}{v_{1}}$,
$g_{1}(t)=k(t)v_{1}^{p-1}(h^{p}-h)/(h-1)$ when $h\neq1$, and
$g_{1}(t)=(p-1)k(t)v_{1}^{p-1}$ when $h=1$, $g_{2}(t)=\mu
f(e^{t})e^{(m+2)t}/v_{1}, g(t)=g_{1}(t)-g_{2}(t)$. Then $Q(t)$
satisfies
$$Q''-f(-t)Q'+g(-t)Q=0.$$ Suppose $Q(t)\not \equiv0$. By Corollary
2.3 and Lemma 3.2, $\lim_{t\to-\infty}Q(t)=0, \lim
_{t\to-\infty}Q'(t)=0, \lim_{t\to+\infty} f(-t)=b_{0}$ and
$\lim_{t\to+\infty}g(-t)=c_{0}$. But by the Lemma 3.3, $Q'(t)$ is
unbounded as $t\to+\infty$. The contradiction shows (1.2) has only
one singular solution.

Now we show the existence of singular solution of (1.2). By Theorem
A, for any $\alpha>\alpha_{**}$, there exists a positive solution
$u_{\alpha}$ of (1.2) which is strictly increasing in $\alpha$, and
$r^{m}u_{\alpha}<\frac{L}{[r^{-l}K(r)]^{\frac{1}{p-1}}}$ by Lemma
2.2 of [6]. From (1.2) we have that
$$(r^{n-1}u_{\alpha}'(r))'+r^{n-1}(K(r)u_{\alpha}^{p}+\mu f(r))=0$$
Integrating from 0 to $r$, we have that, from the fact $r^{-l}K(r)$
is non-increasing,
\begin{eqnarray*}
u'_{\alpha}(r)&=&-\frac{1}{r^{n-1}}\int_{0}^{r}(K(s)u_{\alpha}^{p}(s)+\mu f)s^{n-1}\\
&\leq&\frac{L^{p}}{r^{n-1}}\int_{0}^{r}s^{n-1-\frac{2p}{p-1}}K(s)^{-\frac{1}{p-1}}ds
+\frac{1}{r^{n-1}}(\int_{0}^{r_{0}}+\int_{r_{0}}^{r}\mu
f(s)s^{n-1}ds\\
&\leq&
\frac{L^{p}}{r^{n-1}}r^{\frac{l}{p-1}}K(r)^{\frac{-1}{p-1}}\int_{0}^{r}s^{n-1-\frac{2p}{p-1}-\frac{l}{p-1}}ds
+(Cr^{1-n}+Cr^{1-q})\\
&=&\frac{(p-1)L^{p}}{[(n-2)p-(n+l)][r^{p+1}K(r)]^{\frac{1}{p-1}}}+(Cr^{1-n}+Cr^{1-q}),
\end{eqnarray*}
where $r_{0}$ is positive number. Hence, $u_{\alpha}'$ is uniformly
bounded on any compact subset of $(0,\infty)$ for $\alpha$ and
consequently, ${u_{\alpha}}$ is equicontinuous on any compact subset
of $(0,\infty)$. By the Arzela-Ascoli Theorem and a standard
diagonal argument, there exists a sequence ${\alpha_{j}}$ tending to
$\infty$ as $j\to\infty$ such that $u_{\alpha_{j}}$ converges
uniformly on compact subsets, and thus
$U(r):=\lim_{\alpha_{j}\to\infty}u_{\alpha_{j}}$ is continuous on
$(0,\infty)$. From the Lemma 2.2 of [6], we have that, for each
$\alpha>0$, $$r^{m}u_{\alpha}<r^{m}U(r)\leq
L/(r^{-l}K(r))^{\frac{1}{p-1}}\ \ \mbox{on}\ \ (0,\infty).$$
 Considering the equation
$$u''_{\alpha}=-\frac{n-1}{r}u'_{\alpha}-Ku_{\alpha}^{p}-\mu f,\eqno(3.1)$$ we
observe that $u_{\alpha}''$ is uniformly bounded in $\alpha$ on any
compact subset of $(0,\infty)$. The Arzela-Ascoli Theorem implies
again that there exists a subsequence of ${\alpha_{j}}$ (still
denote it by ${\alpha_{j}}$) such that $u'_{\alpha_{j}}(r)$
converges uniformly  on any compact subset of $(0,\infty)$ as
$\alpha_{j}\to\infty$. Then, $U$ is differentiable on $(0,\infty)$
and $u'_{\alpha_{j}}(r)\to U'$ uniformly on any compact subset of
$(0,\infty)$ as $\alpha_{j}\to\infty$. By (3.1), $u''_{\alpha_{j}}$
converges also uniformly on any compact subset. Therefore, $U'$ is
differentiable on $(0,\infty)$ and $u''_{\alpha_{j}}\to U''$
uniformly on any compact subset of $(0,\infty)$ as
$\alpha_{j}\to\infty$. From (3.1), letting $j\to\infty,U\geq0$
satisfies
$$U''=-\frac{n-1}{r}U'-KU^{p}-\mu f\ \ \ \mbox{on}\ \ (0,\infty)$$ and is a
singular solution. By the maximum principle, $U>0$, which completes
the proof.\vskip0.1in

\textbf{Acknowledgement}\vskip0.1in

The authors would like to express their thanks to Prof Yi-Li for
some helpful suggestions and the referee for helpful comments.

\end{document}